\documentclass[12pt]{amsart}
\usepackage[top=70pt,bottom=70pt,left=75pt,right=77pt]{geometry}

\usepackage{amssymb, amsmath, amsthm, setspace, bbm, prettyref, graphicx, float, subfigure, color, mathrsfs, mathtools,comment,float}
\usepackage{Lutsko_Style}
\DeclareMathSizes{10}{10}{7}{5}

\title{Sarnak's Spectral Gap Question}
\author{Dubi Kelmer}
\thanks{Kelmer is partially supported by NSF CAREER grant DMS-1651563.}
\email{kelmer@bc.edu}
\address{Department of Mathematics, Boston College, Boston, MA}
\author{Alex Kontorovich}
\thanks{Kontorovich is partially supported by NSF grant DMS-1802119 and BSF grant 2020119.}
\email{alex.kontorovich@rutgers.edu}
\address{Department of Mathematics, Rutgers University, New Brunswick, NJ}
\author{Christopher Lutsko}
\email{chris.lutsko@rutgers.edu}
\address{Department of Mathematics, Rutgers University, New Brunswick, NJ}

\newcommand{\Isom}{\operatorname{Isom}}
\newcommand\pf{\begin{proof}}
\newcommand\epf{\end{proof}}

\begin{document}

  \maketitle
\begin{abstract}
\noindent    
We answer in the affirmative a question of Sarnak's from 2007, confirming that
the Patterson-Sullivan base eigenfunction is the unique square-integrable eigenfunction of the hyperbolic Laplacian invariant under the group of symmetries of the Apollonian packing. Thus the latter has a maximal spectral gap. We prove further restrictions on  the spectrum of the Laplacian on a wide class of manifolds coming from Kleinian sphere packings. 
\end{abstract}

  \onehalfspacing
  \setlength{\abovedisplayskip}{1mm}

  \section{Introduction}\label{s:Introduction}

In 2007, Sarnak \cite{Sarnak2007} asked whether the Patterson-Sullivan base eigenvalue (see \S\ref{sec:prelims} for definitions and background) exhausts the discrete spectrum of the hyperbolic Laplacian of the Apollonian 3-fold, so that this manifold has maximal ``spectral gap''; he suggested that this may indeed be the case. The purpose of this paper is to answer this question in the affirmative:

\begin{theorem}[Spectral Gap for the Apollonian Group]\label{thm:Sarnak}
    Let $\Gamma < \Isom(\half^3)$ be the symmetry group of an Apollonian circle packing. Then the base eigenvalue $\lambda_0\approx 0.9065\dots$ is the only discrete eigenvalue of the Laplacian acting on $L^2(\Gamma\bk\half^3)$. 
\end{theorem}

This result implies an improved counting estimate on the number of circles in an Apollonian packing with curvature bounded by a growing parameter (see \cite[Theorem 1]{KontorovichLutsko2022}).

\begin{corollary}
    For a fixed bounded Apollonian circle packing, $\cP$, the number $\cN(T)$ of circles in $\cP$ with curvature bounded by $T$ has the following asymptotic form:
    $$
    \cN(T) = c_\cP\, T^\gd + O (T^\eta (\log T)^{2/5}),
    $$
    as $T\to \infty.$
    Here
    $\gd\approx 1.30\dots$ is the Hausdorff dimension of the residual set of $\cP$, the error exponent is $$\eta=\frac35\delta+\frac25\approx1.18\dots,$$
    and
    $c_\cP>0$ is a constant depending on $\cP$.
\end{corollary}

In general, to each transitive Kleinian packing  (see \S \ref{sec:prelims} for definitions), one can associate a lattice (the supergroup). Theorem \ref{thm:Sarnak} is a consequence of the following more general theorem.

\begin{theorem}\label{thm:general}
    Let $\cP$ be a transitive Kleinian sphere packing  with symmetry group $\Gamma < \Isom(\half^{n+1})$ and associated supergroup $\wt{\Gamma}$. Let $k$ (resp. $\wt{k}$) denote the number of discrete eigenvalues of the Laplacian acting on $L^2(\Gamma \bk \half^{n+1})$ (resp. $L^2(\wt\Gamma \bk \half^{n+1})$) below $n-1$; then $k \le \wt{k}$.
\end{theorem}

Theorem \ref{thm:Sarnak} then follows from combining two facts: $(a)$  that
for circle packings (that is, when $n=2$), $n-1=(n/2)^2=1$, and $(b)$ that in the Apollonian case, $\wt{\Gamma}$ is a double cover of $\PSL(\Z[i])$; the latter satisfies Selberg's eigenvalue conjecture.
The same argument applies to several other packings appearing in the literature to prove maximal spectral gaps.

To illustrate the methodology, we give an alternate proof of the theorem of Phillips-Sarnak \cite{PhillipsSarnak1985} that infinite volume Hecke ``triangle'' groups similarly have only the base eigenvalue in their discrete Laplace spectrum. 

\begin{theorem}[{\cite[Theorem 6.1]{PhillipsSarnak1985}}]\label{thm:PS}
	For $\mu>2$, let  $\Gamma_\mu := \<\mattwos 1\mu01, \mattwos0{-1}10\>\ <\  \Isom(\half^2)$ be an infinite co-volume Hecke triangle group. Then the base eigenvalue $\gl_0(\mu)$ is the unique discrete eigenvalue of the Laplacian acting on $L^2(\Gamma_\mu\bk\half^2)$.
\end{theorem}

One of the key ideas is to replace Dirichlet boundary conditions and nodal domains with Neumann boundary conditions and carefully chosen reflective walls, and to appeal to the fact that the Hecke triangle group with $\mu=2$ (which is a congruence subgroup of the modular group) is known to satisfy Selberg's eigenvalue conjecture. After some preliminaries in \S\ref{sec:prelims}, we give the proofs of these theorems in \S\ref{sec:pfs}.

\subsection*{Acknowledgements}
We thank Peter Sarnak for bringing this problem to our attention, and comments on an earlier draft.

\section{Preliminaries}\label{sec:prelims}

\subsection{Kleinian Packings}

\begin{figure}
\includegraphics[width=3in]{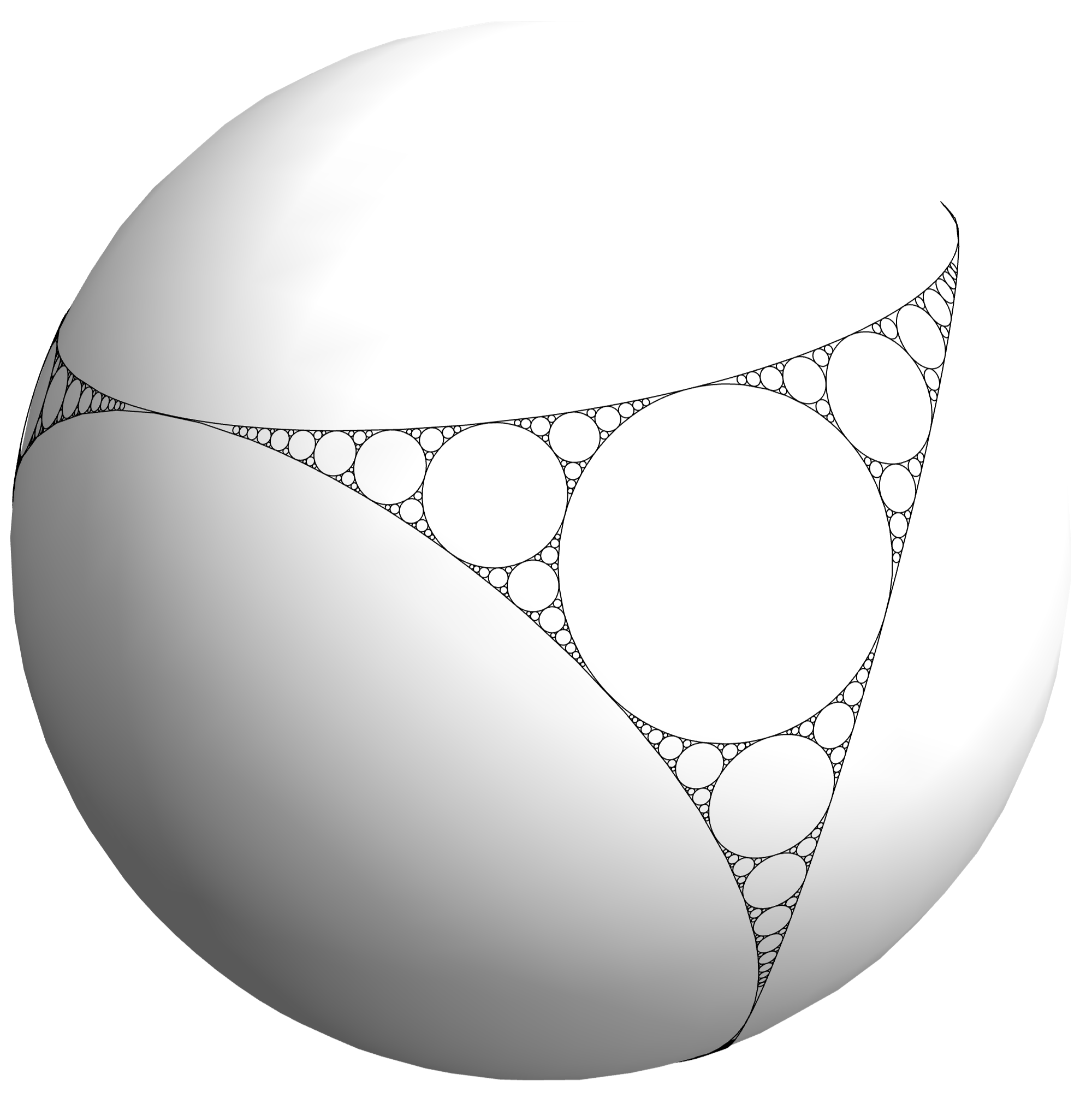}
\caption{An Apollonian circle packing. (Image by Iv\'an Rasskin.)}
\label{fig:apoll}
\end{figure}

The classical Apollonian circle packing is shown in Figure \ref{fig:apoll}. More generally, by a sphere packing $\cP$ of the $n$-sphere, $\bS^n$, we mean an infinite collection of round balls $B$ in $\bS^n$ with pairwise disjoint interiors, such that their closure is all of $\bS^n$. We treat $\bS^n$ as the ideal boundary of the ball model for hyperbolic space $\bH^{n+1}$. A ball $B$ in $\bS^n$ is the  boundary at infinity of a half-space bounded by a  hyperplane $H$ in $\bH^{n+1}$; we denote by $R_B\in \Isom(\bH^{n+1})$ the isometry of $\bH^{n+1}$ corresponding to reflection through $H$. Given any packing $\cP$, we define the \emph{reflection group} 
$$
\G_\cP<\Isom(\bH^{n+1})
$$ of $\cP$, to be the group (infinitely) generated by the reflections $R_B$, for all balls $B$ in $\cP$.

We
are interested in those packings $\cP$ which admit a large group of symmetries. In particular, a packing is defined (in \cite{KapovichKontorovich2021}) to be \emph{Kleinian} if the residual set of $\cP$ agrees with the limit set of some such discrete, geometrically finite group $\G$; the latter is called a \emph{symmetry group} of $\cP$.  When $\G$ is generated by reflections in finitely-many hyperplanes, we say that $\cP$ is \emph{crystallographic}; the study of this latter subclass was initiated in  \cite{KontorovichNakamura2019}. In either setting, there is a Structure Theorem describing how such arise.  Given a Kleinian packing $\cP$ with symmetry group $\G$, we let $\widetilde\G$ denote the so-called \emph{supergroup} of $\cP$, which is the group generated by both the symmetry group $\G$ and the reflection group $\G_\cP$ of $\cP$. 

\begin{theorem}[{Structure Theorem \cite{KontorovichNakamura2019, KapovichKontorovich2021}}]\label{thm:Struct}
%
Suppose that we are given a lattice $\widetilde\G<\Isom(\bH^{n+1})$ containing at least one reflection, with a convex fundamental polyhedron $D$ for the action of $\widetilde \G$, and a minimal set $\widetilde S$ of generators for $\widetilde\G$ whose elements correspond to face-pairing transformations of $D$. Let $T\subset \widetilde S$ be a nonempty subset consisting of some reflections in hyperplanes $H_\alpha$, with $\alpha$ running over some indexing set $A$. The ideal boundary of each $H_\alpha$ is a ball $B_\alpha$; we assume the interiors of these balls are pairwise disjoint, and that the hyperplanes $H_\alpha$ meet the other bounding walls of $D$ either tangentially or orthogonally or not at all. Let $S:=\widetilde S\setminus T$ be the generators of $\widetilde\G$ which are not in $T$, and let $\G_S:=\<S\>$. Then the orbit of the balls $B_\alpha$, $\alpha\in A$, under the group $\G_S$, is a Kleinian packing. (Conversely, every Kleinian packing arises via this construction.)
\end{theorem}

 \begin{figure}
  \begin{center}    
    \includegraphics[width=0.45\textwidth]{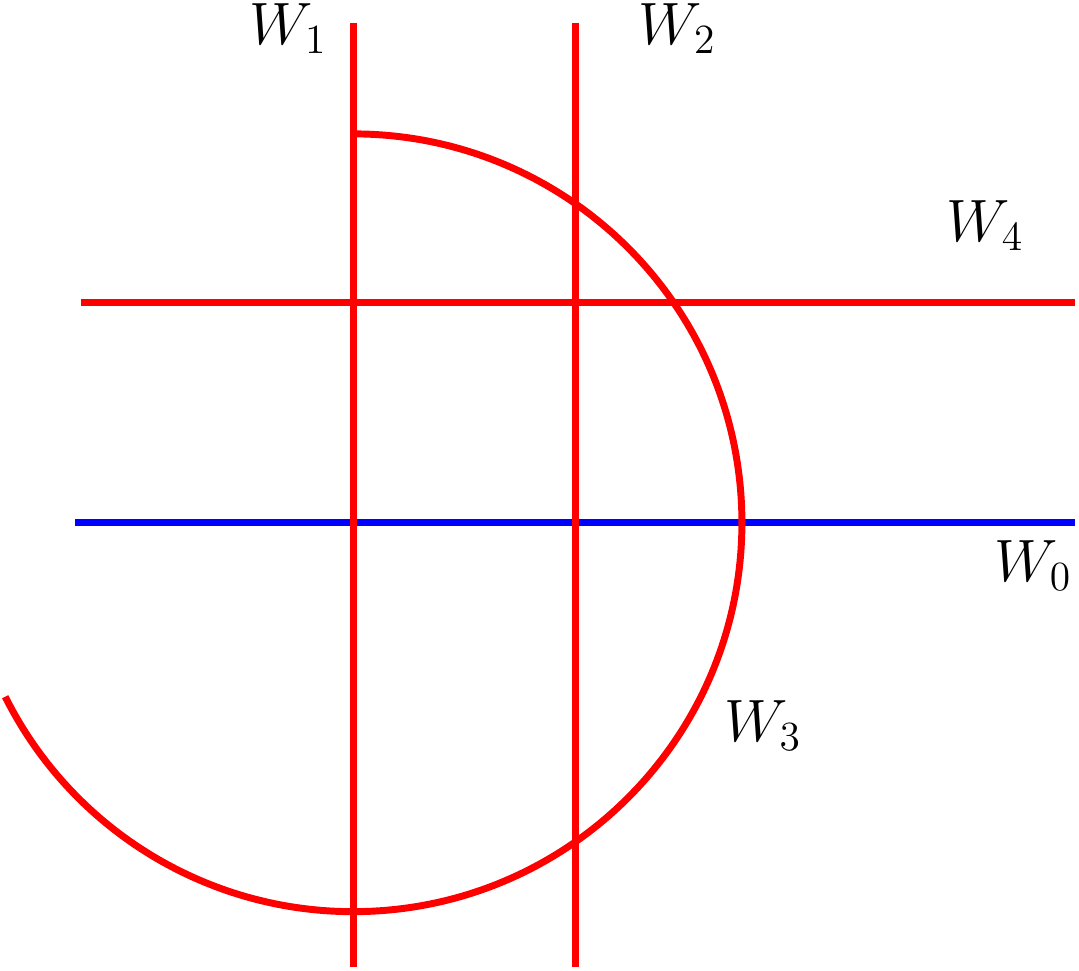}
  \end{center}
  \caption{The Structure Theorem as applied to the Apollonian packing}
  \label{fig:W_0}
\end{figure}
To see how this theorem applies in practice, recall the following construction of the Apollonian packing (which will be relevant to our main theorem). The lines and circles $\{W_\alpha, \alpha=0,\dots,4\}$ illustrated in Figure \ref{fig:W_0} show the boundaries in the plane (identified with $\bS^2$ by stereographic projection and compactification with a point at infinity) of hyperplanes $H_\alpha$; let $\widetilde S$ denote the set of reflections $R_{H_\alpha}$ in these hyperplanes. The group $\widetilde\G:=\<\widetilde S\>$ generated by these reflections is well-known to be a reflective double cover of the Picard group $\PSL_2(\Z[i])$, in the same way that the $(2,3,\infty)$ Hecke triangle group is a reflective double cover of the modular group $\PSL_2(\Z)$. Following the recipe in the Structure Theorem, we set $T$ to be the singleton consisting of inversion in $H_0$; the  boundary of $H_0$, namely, $W_0$, is shown in blue in Figure \ref{fig:W_0}. The set $S:=\widetilde S\setminus T$ of inversions through the other $W_\alpha$ generates the well-known ``Apollonian group'' $\G_S :=\<S\>$, and the orbit of $W_0$ under $\G_S$ is the classical Apollonian packing.

A packing is called {\bf transitive} if its symmetry group takes any sphere to any other (or equivalently, the subset $T\subset\widetilde S$ is a singleton).
The above construction shows that  the Apollonian packing is transitive.

\

\subsection{Spectral Theory}\label{sec:specThPrelim}
Let $n\ge 1$ and let $\half^{n+1}=\{(x_1,\dots,x_n,y): x_j\in\R, \ y>0\}$ denote the upper-half space model of hyperbolic $(n+1)$-space with volume element $\rd V = \frac{\rd x \rd y}{y^{n+1}}$. Further, let $\Delta$ denote the (positive-semidefinite) hyperbolic Laplace-Beltrami operator, which in these coordinates is given by
  \begin{align*}
     \Delta = - y^2 \left( \frac{\partial^2}{\partial y^2} + \frac{\partial^2}{\partial x_1^2} + \dots +  \frac{\partial^2}{\partial x_{n}^2}\right) + (n-1) y \frac{\partial}{\partial y}.
  \end{align*}
Let $\Gamma$ be a discrete, geometrically finite group of isometries of $\half^{n+1}$, let $\Lambda_\Gamma\subseteq \partial\half^{n+1}$ denote the limit set of $\Gamma$, that is, the set of limit points in the ideal boundary of a $\Gamma$-orbit, and let $\delta=\delta_\Gamma$ be the Hausdorff dimension of $\Lambda_\Gamma$. We are interested in the spectrum of the Laplacian acting on the space $L^2(\Gamma\bk \half^{n+1})$ of square-integrable $\Gamma$-invariant functions.

When $\Gamma$ is a non-lattice, that is, it acts on $\half^{n+1}$ with \emph{infinite} co-volume, the work of Lax-Phillips \cite{LaxPhillips1982} shows that there are at most finitely many discrete eigenvalues of $\Delta$ below $(n/2)^2$, and that the spectrum above $(n/2)^2$ is purely continuous. The bottom discrete eigenvalue is given by the Patterson-Sullivan formula \cite{Patterson1976, Sullivan1984}:
$$
\lambda_0 = \delta(n-\delta),
$$
and only exists if $\delta>n/2$; otherwise, the discrete spectrum is empty. 
In general, if we are given a domain $D\subset\half^{n+1}$ so that the Laplace spectrum (possibly with some specified boundary conditions) below $(n/2)^2$ is empty, we say that $D$ is \emph{eigenvalue-free}. For example, a result of Phillips-Sarnak \cite[Theorem 3.7]{PhillipsSarnak1985} states that if $D$ is the region bounded by exactly $m$ geodesic hyperplanes, and 
\be\label{eq:sidesBnd}
m\le \lfloor (n+4)/2 \rfloor,
\ee
then $D$ is eigenvalue-free for the Laplacian with Neumann boundary, that is, having vanishing normal derivatives across the boundary of $D$. (Here the floor function returns the integer part of its argument.)

When $\delta_\Gamma>n/2$, so that the base eigenvalue $\lambda_0$ exists, let $\lambda_1$ denote either the next smallest eigenvalue above the base, or $(n/2)^2$ if none such exist; we call the difference $\lambda_1 - \lambda_0$ the \emph{spectral gap} for $\Gamma$. When $\lambda_1=(n/2)^2$, we say the spectral gap is \emph{maximal}.
As a point of reference, note that the Selberg Eigenvalue Conjecture states, in the context of $\Gamma$ being a (finite-volume) congruence subgroup of the modular group (for which $n=1$, $\delta_\Gamma=1$, and $\gl_0=0$, corresponding to the constant function), that the spectral gap is maximal, that is, it is equal to $1/4$.
This maximal gap is known to hold for many congruence groups \cite{Huxley1996, BookerStrombergssonVenkatesh2006} in $\bH^2$, and also for a number of Biachni groups acting on $\half^3$.
By the latter, we mean groups $\PSL_2(\cO_d)$, where, for $d>0$, $\cO_d$ is the ring of integers of the imaginary quadratic field $\Q(\sqrt{-d})$.
For example, in \cite[Prop. 6.2, p. 407]{ElstrodtGrunewaldMennicke1998} it is shown that:
\be\label{eq:EGM}
\lambda_1(\PSL_2(\cO_d))\ge1,
\ee
for $d=1, 2, 3$, and $7$; that is, these groups all have a maximal gap.

\

\subsection{Calculus of Variations}
Our last  preliminaries involve recording some well-known results from the calculus of variations  applied to eigenvalue problems in the setting of domains in hyperbolic space with its Laplacian. 
Let $D$ be an open connected subdomain of $\bH^{n+1}$ with ``nice'' boundary, and let $W^1(D):=\{f\in L^2(D) \ : \ \nabla f \in L^2(D)\}$.
Suppose that a hyperplane $H$ splits a domain $D$ into two pieces 
$$
D = D_1\sqcup D_2 \sqcup (H\cap D),
$$
with $D_j$ both open, connected, and non-empty. 
Let $\lambda_i(D)$, $i=0,1,\dots,k$, denote the eigenvalues of the hyperbolic Laplacian on $D$ with some specified (say, Neumann or Dirichlet) conditions on the boundary.
Similarly, let $\lambda_i(D_j)$, $i=0,1,\dots,k_j$, denote the eigenvalues on $D_j$, $j=1,2$ with Neumann boundary conditions along the common boundary of $D_j$ with $H$. 
  Then a simple application of the min-max principle yields the following inequality, see, e.g., {\cite[p. 409, Theorem 4]{CourantHilbert1953}}.
  \begin{theorem}\label{thm:varCalc}
   For any $\kappa>0,$
  \begin{align}
    \#\{i \le k  : \lambda_i(D) < \kappa \} \le     \#\{i \le k_1 : \lambda_i(D_1) < \kappa \} + \#\{i \le k_2  : \lambda_i(D_2) < \kappa \}.
  \end{align}
  \end{theorem}
That is, the number of eigenvalues on $D$ below a parameter $\kappa$ is at most the sum of the number of eigenvalues below $\kappa$ on the two regions remaining after a Neumann cut. 
Note further that the Neumann condition arises naturally in the context of reflections in hyperplanes: an eigenfunction invariant under such a reflection must have vanishing normal derivate across the hyperplane.

\section{Proofs of Theorems  \ref{thm:general} and \ref{thm:PS}}\label{sec:pfs}

\subsection{Hecke Groups}

The original approach of Phillips-Sarnak to Theorem \ref{thm:PS} is to study the topology of a hypothetical eigenfunction's nodal domain, along the boundary of which the eigenfunction of course satisfies the Dirichlet condition (that is, vanishes); see \cite[Fig 6.3]{PhillipsSarnak1985}.
This approach seems to become rather cumbersome in the Apollonian setting, as the nodal domain can interact with a fundamental domain for the group action in rather complicated ways. But 
Theorem \ref{thm:varCalc}, together with some carefully chosen additional reflective hyperplanes, solves both problems with ease.

\pf[Proof of Theorem \ref{thm:PS}]

\begin{figure}
  \begin{center}    
    \includegraphics[width=0.4\textwidth]{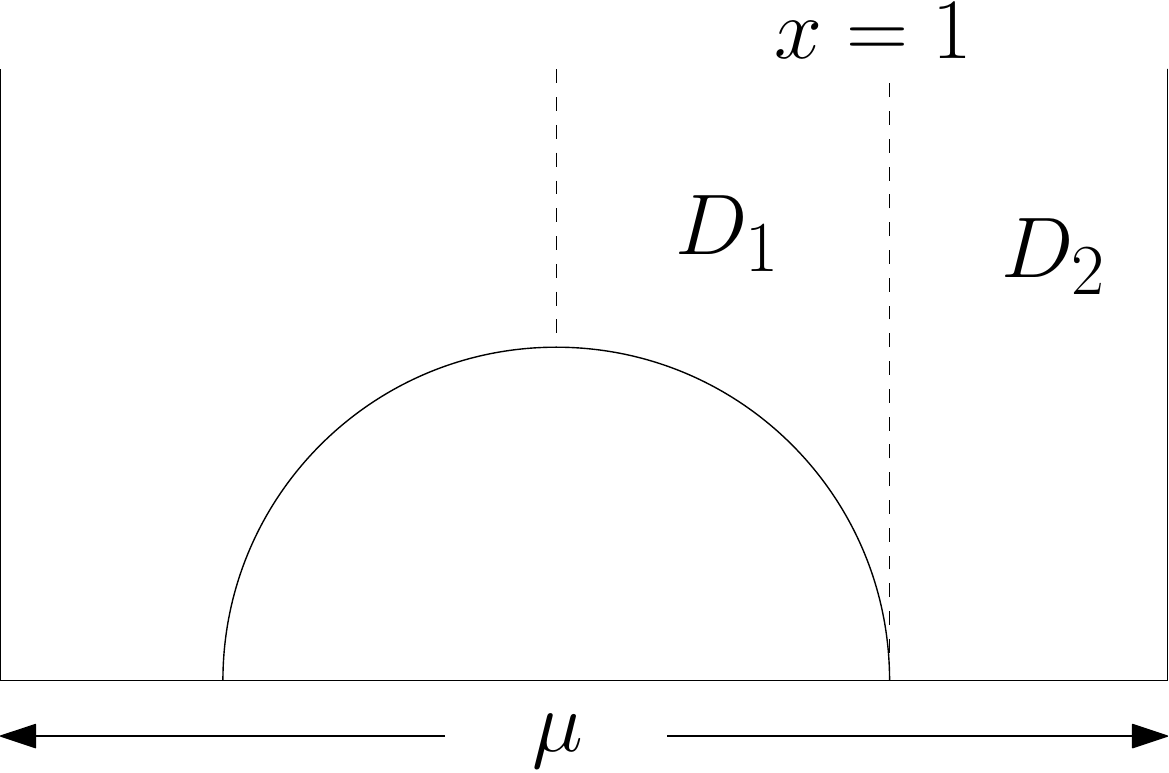}
  \end{center}
  \caption{The regions $D_1$ and $D_2$ in the proof of Theorem \ref{thm:PS}}
  \label{fig:Hecke}
\end{figure}

Fix $\mu>2$ and let $\G_\mu$ be a Hecke triangle group as in the statement of the theorem.
Consider the group $\G<\Isom(\bH^2)$ generated by reflection across the $y$-axis, the unit circle, and the line $x=\mu/2$. Then $\G_\mu$ is  the orientation-preserving (index-two) subgroup of $\G$; in particular, if the latter has no eigenvalues besides the Patterson-Sullivan, the same holds for the former. Let $D$ denote the domain bounded by these three reflections generating $\G$. We cut $D$ with the hyperline $x=1$, splitting it into two domains: $D_1$ bounded by the $y$-axis, unit circle, and line $x=1$, and $D_2$ bounded only by $x=1$ and $x=\mu/2$; see Figure \ref{fig:Hecke}.

We apply Theorem \ref{thm:varCalc} with $\kappa=1/4.$ The region $D_1$ with Neumann boundary is the fundamental domain of a congruence lattice, namely, it contains the level-2 principal congruence subgroup of the modular group. As above \eqref{eq:EGM}, this is well-known to satisfy Selberg's Eigenvalue Conjecture, so the number of eigenvalues below $1/4$ is exactly one, corresponding to the constant function. On the other hand, $D_2$ is bounded by exactly two hyperplanes; by \eqref{eq:sidesBnd}, $D_2$ is eigenvalue-free. Therefore Theorem \ref{thm:varCalc} tells us that $D$ itself has at most one eigenvalue below $1/4$, which we know to be the Patterson-Sullivan base eigenvalue.
\epf

\subsection{Proof of Theorem \ref{thm:general}}

We follow a similar strategy for
sphere packings.
Let $\cP$ be a transitive Kleinian sphere packing with symmetry group $\Gamma$.
Recall from the Structure Theorem \ref{thm:Struct} that $\Gamma$ is obtained from its supergroup $\wt{\Gamma}$ by removing from a minimal set of generators the  reflection through a single hyperplane $H_0$ (see Figure \ref{fig:W_0} for the Apollonian example). Moreover we have a fundamental domain $\wt\cF$ for $\wt\Gamma$ so that $H_0$ only intersects the other walls of $\wt\cF$ at right angles.

 \begin{figure}
  \begin{center}    
    \includegraphics[width=0.4\textwidth]{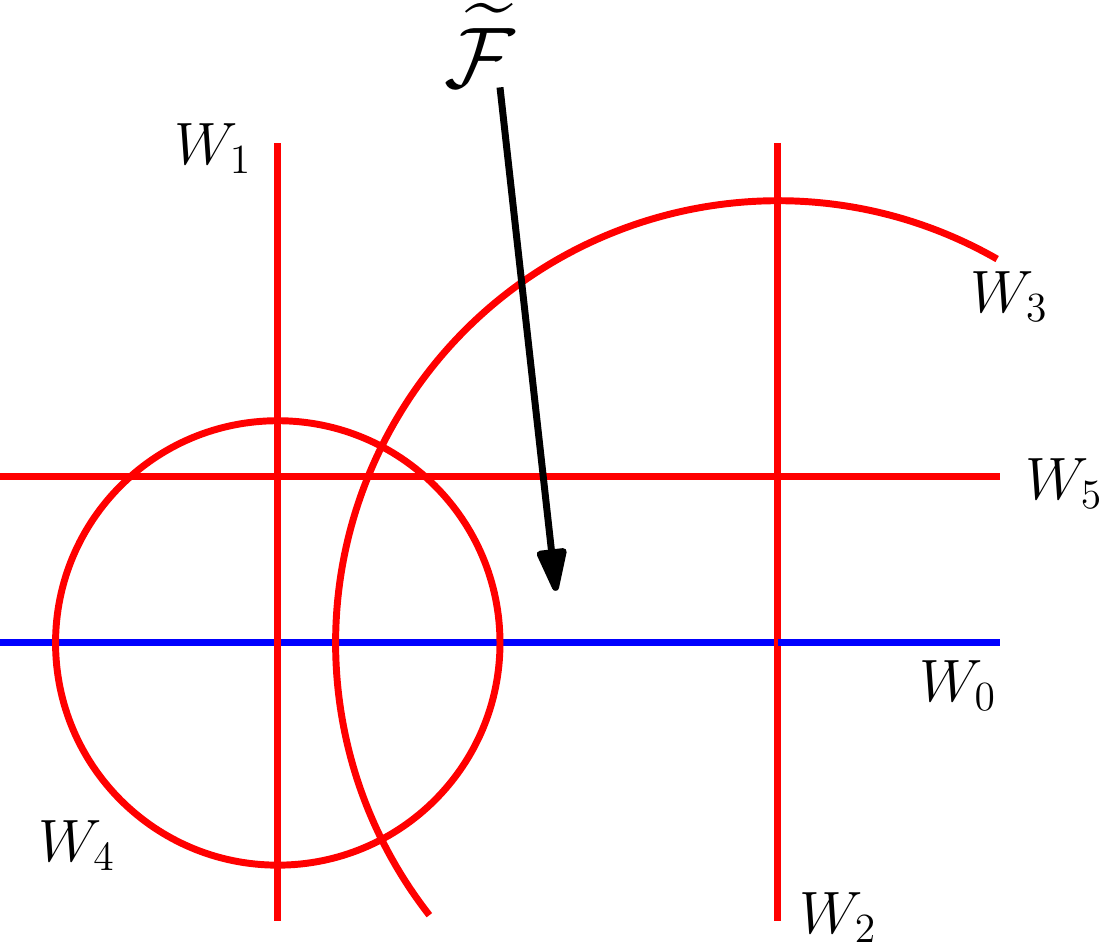}
    \includegraphics[width=0.4\textwidth]{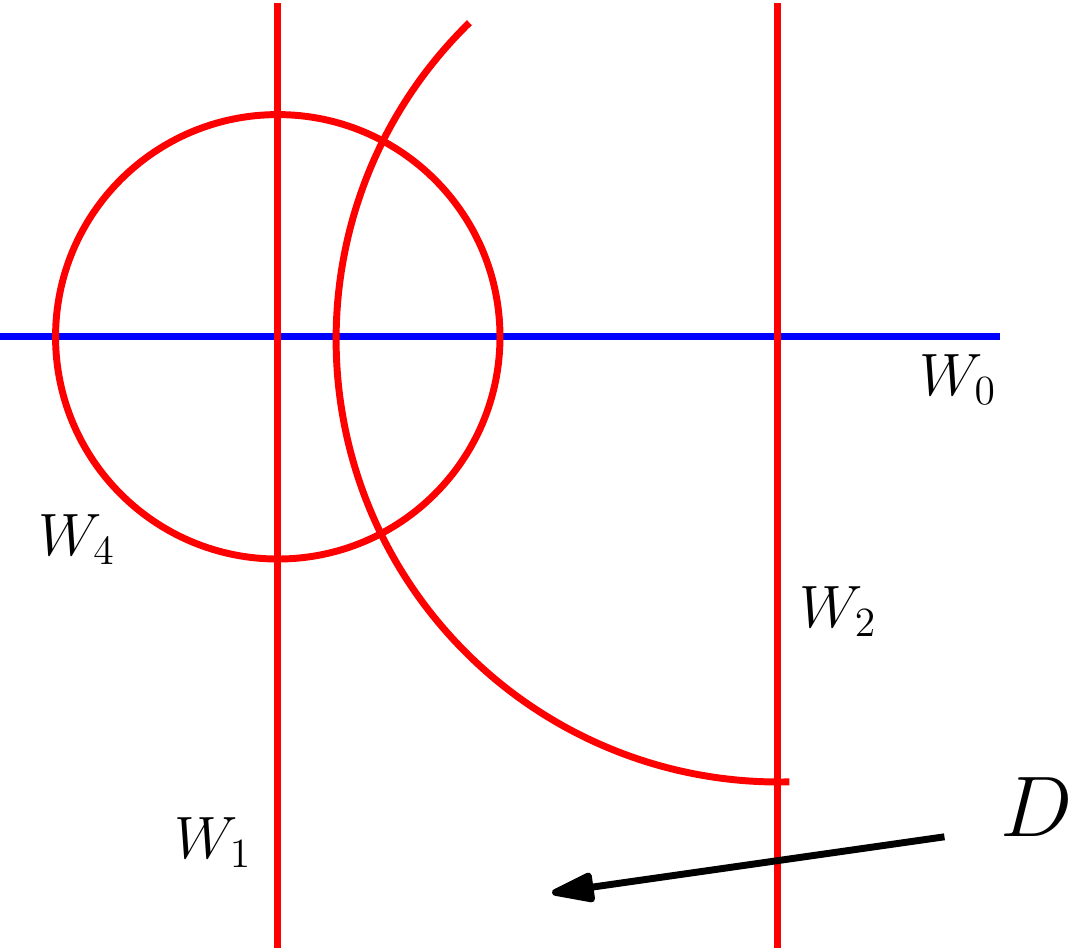}    
  \end{center}
  \caption{The regions $\wt\cF$ and $D$ in the proof of Theorem \ref{thm:general}}
  \label{fig:Ddef}
\end{figure}

    Let $\cF$ denote the corresponding fundamental domain for $\Gamma$, so that  $H_0$ splits $\cF$ into the fundamental domain $\wt{\cF}$ for $\wt{\Gamma}$, and a region $D$, say; see Figure \ref{fig:Ddef} for an example. Thus using Theorem \ref{thm:varCalc} with $\kappa = n-1$, we can upper-bound the number of eigenvalues of $\Gamma$ below $n-1$ by 
    \begin{align}
        \#\{i : \lambda_i(\cF) < n-1\} \le 
        \#\{i : \lambda_i(\wt{\cF}) < n-1\}+
        \#\{i : \lambda_i( D) < n-1\}.
    \end{align}
    Thus it remains to show that $\#\{i : \lambda_i( D) < n-1\}=0$.

For this, notice that, by construction, $D$ is a fundamental domain for the action of the stabilizer, $\Gamma_0$, say, of $H_0$ in $\Gamma$. 
But by the Structure Theorem, $\Gamma_0$ preserves $H_0$ and acts on it as a discrete lattice; hence the limit set of $\Gamma_0$ must be equal to the ideal boundary, $W_0$, of $H_0$. Therefore the Hausdorff dimension $\gd_{\Gamma_0}$ of this limit set is equal to the dimension of $W_0$, which is $n-1$. Applying the Patterson-Sullivan formula shows that $D$ has no eigenvalues below $\gd_{\Gamma_0}(n-\gd_{\Gamma_0})=n-1$, as claimed.
\qed

  \small 
  \bibliographystyle{alpha}
  \bibliography{biblio}

\end{document}